\documentclass[a4paper]{article}

\usepackage{epsf}
\usepackage{amssymb}
\usepackage{amsfonts}
\usepackage{xypic}
\usepackage{headings}

\newcommand{\inter}[1]{\mathord{\mathop{#1}\limits^{\smash{\circ}}}}
\newcommand{\interaccent}[1]{\mathord
                              {\mathop{#1}\limits^{\smash{\circ}}\negthinspace{}'
                            }
                              }

\newcommand{\mono}{\hookrightarrow}
\newcommand{\U}{\tilde{U}}
\newcommand{\id}{\mbox{\rm id}}
\newcommand{\Gl}{\mbox{\rm Gl}}
\newcommand{\St}{\mbox{\rm St}}
\newcommand{\two}{\rightrightarrows}
\newcommand{\Hom}{\mbox{\rm Hom}}

\newtheorem{lemma}{Lemma}[subsection]
\newtheorem{prop}[lemma]{Proposition}
\newtheorem{stel}[lemma]{Theorem}
\newtheorem{gevolg}[lemma]{Corollary}

\newcounter{lijst}

\newenvironment{mylist}{
\setcounter{lijst}{1}
\begin{list}{
             \parbox{4ex}{(\roman{lijst})}
} 
{\usecounter{lijst}\itemsep=0 pt \topsep=0 pt \leftmargin=0 pt\leftmargini=0pt \parsep=0pt \parskip=0pt\labelwidth=0pt \labelsep=0pt \listparindent=3ex
}}{\end{list}}

\begin{document}
\title{Simplicial Cohomology of Orbifolds}
\author{I. Moerdijk and D.A. Pronk\thanks{This research is part of a project funded by the Netherlands Science Organization (NWO).}}
\renewcommand{\today}{}
\maketitle

\begin{abstract} For any orbifold ${\cal M}$, we explicitly construct a
simplicial complex
$S({\cal M})$ from a given triangulation
of the `coarse' underlying space together with
the local isotropy groups of ${\cal M}$. We prove that, for any local
system on ${\cal M}$,
this complex $S({\cal M})$ has the 
same cohomology as ${\cal M}$. The use of $S({\cal M})$ in explicit
calculations is illustrated in the example of the `teardrop'
orbifold.
\end{abstract}

\section*{Introduction.}
Orbifolds or V-manifolds were first introduced by Satake \cite{sa56}, 
and arise
naturally in many ways. For example, the orbit space of any proper action by a
(discrete) group on a manifold has the structure of an orbifold; this applies
in particular to moduli spaces. Furthermore, the orbit space of any almost free
action by a compact Lie group has the structure of an orbifold, as does the
leaf space of any foliation with compact leaves and finite holonomy. 
Examples of orbifolds are discussed in \cite{DM,sa56,Thurston} and many others.

For an orbifold ${\cal M}$, one can define in a natural way a cohomology theory 
with
coefficients in any local system on ${\cal M}$. This cohomology is not an 
invariant of
the underlying (`coarse') space, but of the finer orbifold structure. If the
orbifold is given as the orbit space $X/G$ of a group action as above, this
cohomology is the equivariant sheaf cohomology of the group action. It agrees
with the (ordinary) cohomology of the Borel construction $EG\times_G X$.

This cohomology is the most natural one for orbifolds. It fits in well with the
notion of fundamental group described in \cite{Thurston}, by the familiar 
`Hurewicz formula' $H^1({\cal M},A) = \Hom(\pi_1({\cal M}),A)$ (where $A$ 
is any abelian group).

The purpose of this paper is to give a {\it simplicial} description of these
cohomology groups, suitable for calculations. More precisely, using
triangulations of singular spaces \cite{gor}, we will associate to any orbifold
${\cal M}$, presented by an orbifold atlas as in \cite{sa56}, a simplicial set $S({\cal M})$. The
construction of $S({\cal M})$ uses the simplices in a triangulation of the 
coarse
underlying space $M$ of ${\cal M}$, as well as all the local isotropy groups. 
The
construction will have the following property.

\paragraph{}{\bf Theorem} {\it For any local system of coefficients ${\cal A}$ 
on the orbifold ${\cal M}$, there is a
canonically associated local system $A$ on the simplicial set $S({\cal M})$, 
for which
there is a natural isomorphism $H^*({\cal M},{\cal A}) \cong H^*(S({\cal M}),A)$.}

\paragraph{}After reviewing some preliminary definitions, we will present 
our construction
in Section~\ref{simpcom} of this paper. The proof of the theorem will be 
based on the fact
that the simplicial set $S({\cal M})$ associated to an orbifold ${\cal M}$ 
is also closely
related to the representation of ${\cal M}$ by a groupoid $G({\cal M})$ 
suggested in
\cite{hae}. In fact, our proof shows that $S({\cal M})$ has the same 
homotopy type as
the classifying space of $G({\cal M})$.

Before giving the proof in Section~\ref{seven}, we will present an
example of calculations based on the simplicial construction for the `teardrop' orbifold. We believe that
much more work should be done in this direction. In fact, the explicit
description of the simplicial set in terms of an atlas for the orbifold, and
the resulting description of the cohomology groups by generators and relations,
makes it suitable for computer assisted calculations.

\section{Preliminaries.}\label{prel}

\subsection{Basic definitions.}\label{one}
In this section we briefly review the basic definitions concerning orbifolds, 
or V-manifolds in the terminology of Satake (see \cite{sa56, sa57, Thurston}).
Let $M$ be a paracompact Hausdorff space. An {\it orbifold chart} on $M$ is 
given by a connected open subset $\tilde{U}\subseteq{\mathbb R}^n$ for some
integer $n\geq 0$,  a finite group $G$ of $C^\infty$-automorphisms of $\tilde{U}$, 
and a map $\varphi\colon\tilde{U}\rightarrow M$, such that $\varphi$ is 
$G$-invariant ($\varphi\circ g=\varphi$ for all $g\in G$) and induces a 
homeomorphism of $\tilde{U}/G$ onto the open subset 
$U=\varphi(\tilde{U})\subseteq M$. An {\it embedding} $\lambda\colon(\tilde{U},G,\varphi)\mono(\tilde{V},H,\psi)$ between two such 
charts is a smooth embedding $\varphi\colon\tilde{U}\mono\tilde{V}$ 
with $\psi\circ\lambda=\varphi$. An {\it orbifold atlas} on $M$ is a 
family ${\cal U}=\{(\U,G,\varphi)\}$ of such charts, which cover $M$ and are 
locally compatible in the following sense: given any two charts $(\U,G,\varphi)$ 
for
$U=\varphi(\U)\subseteq M$ and $(\tilde{V},H,\psi)$ for $V\subseteq M$, and a 
point $x\in U\cap V$, there exists an open neighborhood $W\subseteq U\cap V$ of 
$x$ and a chart $(\tilde{W},K,\chi)$ for $W$ such that there are embeddings
$(\tilde{W},K,\chi)\mono(\U,G,\varphi)$ and $(\tilde{W},K,\chi)\mono(\tilde{V},H,\psi)$. Two such atlases are said to 
be equivalent if they have a common refinement. An {\it orbifold} 
(of dimension $n$) is such a space $M$ with an equivalence class of atlases 
${\cal U}$.
We will generally write ${\cal M}=(M,{\cal U})$ for the orbifold 
${\cal M}$ represented by the space $M$ and a chosen atlas ${\cal U}$.

\subsubsection{Remarks.}\label{rems}\begin{mylist}
\item \label{two}
For two embeddings $\lambda,\mu\colon(\U,G,\varphi)\two(\tilde{V},H,\psi)$ 
between charts, there exists a unique $h\in H$ such that $\mu=h\circ\lambda$. 
In particular, since each $g\in G$ can be viewed as an embedding of $(\tilde{U},G,\varphi)$ into itself, there exists for the two embeddings 
$\lambda$ and $\lambda\circ g$ a unique $h\in H$ with 
$\lambda\circ g=h\circ \lambda$. This $h$ will be denoted by $\lambda(g)$. 
In this way, every embedding $\lambda$ also induces an injective group 
homomorphism, (again denoted) $\lambda\colon G\rightarrow H$, with defining 
equation 
$$
\lambda(g\cdot\tilde{x})=\lambda(g)\lambda(\tilde x)\qquad\qquad (\tilde{x}\in\tilde{U}).
$$
Furthermore, if $h\in H$ is such that $\lambda(\U)\cap h\cdot\lambda(\U)\neq\emptyset$, then $h$ belongs to the image of this 
group homomorphism $\lambda\colon G\rightarrow H$, and hence $\lambda(\tilde{U})=h\cdot\lambda(\U)$.
(This is proved in \cite{sa56} for the codimension 2 case, and in \cite{MP} for 
the general case.)
\item
By the differentiable slice theorem for smooth group actions \cite{ko}, any 
orbifold of dimension $n$ has an atlas consisting of `linear' charts, i.e.\ 
charts of the form $(\U,G)$ where $G$ is a finite group of linear transformations 
and $\U$ is an open ball in ${\mathbb R}^n$.
\item
If $(\U,G,\varphi)$ and $(\tilde{V},H,\psi)$ are two charts for the same 
orbifold structure on $M$, and $\U$ is simply connected, then there exists an
embedding $(\U,G,\varphi)\mono(\tilde{V},H,\psi)$ whenever $U\subseteq V$, 
i.e.\ when $\varphi(\U)\subseteq\psi(\tilde{V})$ (see \cite{sa57}, footnote 2).
{\it In this paper we will take all charts to be simply connected.}
\item\label{contract}
If $(\U,G,\varphi)$ is a chart, and $V$ is a connected (and simply connected)
open subset of $U\subseteq M$, then $V$ inherits a chart structure from $U$ in 
the following way: let $\tilde{V}$ be a connected component 
of $\varphi^{-1}(V)\subseteq\tilde{U}$, and let 
$H=\{g\in G| \,g\cdot\tilde{V}=\tilde{V}\}$. Then $(\tilde{V},H,\varphi|_{\tilde{V}})$
is a chart, which embeds into $(\tilde{U},G,\varphi)$, and hence defines the 
same orbifold structure on points in $V$.
\end{mylist}

\subsubsection{Examples.} We list some well-known examples, see 
e.g.\ 
\cite{hae, MP, Thurston}.
\begin{mylist}
\item
If a discrete group $\Gamma$ acts smoothly and properly on a manifold $N$, 
the orbit space $M=N/\Gamma$ has a natural orbifold structure. Examples 
are
weighted projective space and moduli spaces.
\item 
If a compact Lie group $L$ acts smoothly on a manifold $N$ and the isotropy 
group $L_x$ at each point $x\in N$ is finite and acts faithfully on a slice 
$S_x$ through $x$, then the orbit space $N/L$ has a natural orbifold 
structure.
Moreover any orbifold can be represented this way.
\item
If $N$ is a manifold equipped with a foliation ${\cal F}$ of codimension $n$, 
with the property that each leaf is compact and the holonomy group at each 
point
is finite, then the space of leaves $N/{\cal F}$ has the natural structure
of an $n$-dimensional orbifold. Again, any orbifold can be represented in this way.
\end{mylist}

\subsection{Triangulation of orbifolds.}\label{trian}
Let ${\cal M}=(M,{\cal U})$ be an orbifold of dimension $n$. For each point 
$x\in M$, one can choose a linear chart $(\tilde{U},G,\varphi)$ around $x$, 
with $G$ a finite subgroup of the linear group $\Gl(n,{\mathbb R})$. 
Let $\tilde{x}\in \U$ be a point with $\varphi(\tilde{x})=x$, and 
$G_{x}=\{g\in G|\,g\cdot\tilde{x}=\tilde{x}\}$ the isotropy subgroup at 
$\tilde{x}$. Up to conjugation, $G_x$ is a well defined subgroup of 
$\Gl(n,{\mathbb R})$. The space $M$ thus carries a well-known natural 
stratification, whose strata are the connected components of the sets
$$S_H=\{x\in M|\,(G_x)=(H)\},$$
where $H$ is any finite subgroup of $\Gl(n,{\mathbb R})$ and $(H)$ is its  
conjugacy class.
It is also well-known (\cite{Yang, gor}) that there exists a triangulation 
${\cal T}$ of $M$ subordinate
to this stratification (i.e. the closures of the strata lie on subcomplexes 
of the triangulation ${\cal T}$). By replacing ${\cal T}$ by a stellar 
subdivision, one can assume that the cover of closed simplices in ${\cal T}$ 
refines the atlas ${\cal U}$.
For a simplex in such a triangulation, the isotropy groups of all interior 
points of $\sigma$ are the same, and are subgroups of the isotropy groups of 
the boundary points. By taking a further subdivision of
${\cal T}$, we may in fact assume that there is one face $\sigma'\subset\sigma$ 
such that the isotropy is constant on $\sigma-\sigma'$, and possibly 
larger on
$\sigma'$. In particular, any simplex $\sigma$ will then have a vertex 
$v\in\sigma$ with maximal isotropy,
i.e.\ $G_x\subseteq G_v$ for all $x\in\sigma$. 
We call such a triangulation {\it adapted to} ${\cal U}$. For reference we state:

\begin{prop}
For any orbifold $M$ and any orbifold atlas ${\cal U}$ there exists an 
adapted triangulation ${\cal T}$ for ${\cal M}=(M,{\cal U})$.
\end{prop}

Let $\sigma$ be a (closed) $n$-simplex in such a triangulation,
which is contained in a chart $U$. (In this paper we will take all simplices 
to be closed.)
The next lemma describes how to
lift $\sigma$ to a simplex $\tilde{\sigma}$ in $\U$ by choosing
a connected component of the inverse image $\varphi^{-1}(\inter\sigma)$ of 
the interior.

\begin{lemma}\label{lift}
For every connected component $S$ of the inverse image
$\varphi^{-1}(\inter{\sigma})$ of the interior of $\sigma$, the map 
$\varphi$ restricts to a homeomorphism on the closures $\varphi|_{\overline{S}}\,\colon\overline{S}\stackrel{\sim}%
{\longrightarrow}\overline{\sigma}=\sigma$; in particular the 
triangulation of ${\sigma}$ lifts to a triangulation of $\overline{S}$.
\end{lemma}

\paragraph{Proof.}
As is well-known (see e.g. \cite[Lemma 1]{Yang}), the map 
$\varphi\colon\U\rightarrow\U/G=U$ restricts to a
covering projection on each stratum. In particular, since the isotropy 
is constant on the interior $\inter{\sigma}$, we find that 
$\varphi^{-1}(\inter{\sigma})=\coprod S_i$ is a disjoint sum of open 
simplices $S_i$ with 
$\varphi\colon S_i\stackrel{\sim}{\longrightarrow}\inter{\sigma}$.
Let $S$ be one of these $S_i$.
By continuity, $\varphi$ maps $\overline{S}$ to
${\sigma}$. This restriction is also surjective, because 
${\sigma}\subseteq U$
and the action of $G$ on $\U$ is continuous. 
Since ${\cal T}$ is subordinate to the stratification,
the isotropy groups of the boundary points of $\sigma$
contain the isotropy group of the interior points,
and therefore $\varphi$ has to be one-one on $\overline{S}$. Since 
$\varphi$ is also open (as quotient map), it follows that $\varphi|_{\overline{S}}\,\colon\overline{S}\stackrel{\sim}{\longrightarrow}{\sigma}$.

\paragraph{}By taking further stellar subdivisions of ${\cal T}$, 
we may assume that for every simplex $\sigma$ in ${\cal T}$, the closure of 
the open star
$$
\St (\sigma)=\cup\{\inter{\tau}|\sigma\subseteq\tau\}
$$
is contained in a chart of the atlas ${\cal U}$.
We will call a triangulation of ${\cal M}$ with these properties a {\it good} 
triangulation for ${\cal M}=(M,{\cal U})$.

\begin{prop}
For any orbifold $M$ and any orbifold atlas ${\cal U}$ there exists a good 
triangulation ${\cal T}$ for ${\cal M}=(M,{\cal U})$.
\end{prop}

Consider such a good triangulation. Let $x$ be a point of $M$, let $\sigma(x)$ 
be the smallest simplex containing $x$, and consider the open star neighborhood
$\Sigma_x=\St(\sigma(x))$. Let $(\U,G,\varphi)$ be a chart in the atlas 
${\cal U}$ for which $\overline{\St(\sigma(x))}\subset U=\varphi(\U)$. For later
purposes, it is useful to describe how $\overline{\St(\sigma(x))}$
lifts to a triangulation in $\U$.

\begin{lemma}\label{stars}
The closed star $\overline{\St(\sigma(x))}$ lifts to a triangulation ${\cal R}$
of 
$$
\varphi^{-1}(\overline{\St(\sigma(x))})\subseteq\U
$$
with the property that every connected component of $\varphi^{-1}({\St(\sigma(x))})$ 
contains precisely one lifting $\tilde{x}$ of $x$, and these components are 
open star neighborhoods $\Sigma_{\tilde{x}}=\St(\sigma(\tilde{x}))$ for the 
triangulation ${\cal R}$:
$$
\varphi^{-1}(\Sigma_x)=\coprod_{\tilde{x}\in\varphi^{-1}(x)}\Sigma_{\tilde{x}}.
$$ 
\end{lemma}

\paragraph{Proof.}
According to Lemma~\ref{lift}, for every $n$-simplex 
$\tau\in\overline{\St(\sigma(x))}$ and every connected component of 
$\varphi^{-1}(\inter\tau)$ there is precisely one lifting $\tilde\tau$ of 
$\tau$ containing that connected component.
Since we could lift the whole triangulation of $\tau$ to $\tilde\tau$, 
we get together with $\tilde\tau$ also
all its boundary simplices.
Let ${\cal R}$ be the set of all those liftings (for all $n$-simplices 
in $\overline{\St(\sigma(x))}$ and their boundaries). It is clear that 
they cover $\varphi^{-1}(\overline{\St(\sigma(x))})$. Remark that for 
every
simplex $\rho\in{\cal R}$, $\varphi|\rho$ becomes a homeomorphism
onto a simplex in $\overline{\St(\sigma(x))}$.
Let $\rho_1,\rho_2\in{\cal R}$ be two simplices with 
$\rho_1\cap\rho_2\neq\emptyset$.  Then $\varphi(\rho_1\cap\rho_2)\subseteq\varphi(\rho_1)\cap\varphi(\rho_2)$.
This last intersection is a simplex in ${\cal T}$, which we will denote by $\tau$. 
Let
$$
\tilde{\tau}_i:=\varphi^{-1}(\varphi(\rho_1)\cap\varphi(\rho_2))\cap\rho_i
$$
be the two liftings of $\tau$ in the $\rho_i$. 
(Note that $\rho_1\cap\rho_2=\tilde\tau_1\cap\tilde\tau_2$.) There exists an 
element $g\in G$ such that $g\cdot\tilde\tau_1=\tilde\tau_2$.
This element $g$ is in the isotropy group of $\varphi(\rho_1\cap\rho_2)$, but
not in the isotropy group of the rest of $\tau$. Since ${\cal T}$ is subordinate 
to the stratification induced by the isotropy groups, it follows that
$\varphi(\rho_1\cap\rho_2)$ is a subsimplex of $\tau$, and therefore a subsimplex 
of $\varphi(\rho_1)$.
We know already that 
$\varphi|\rho_1\colon\rho_1\stackrel{\sim}{\longrightarrow}\varphi(\rho_1)$, so $\rho_1\cap\rho_2=\varphi^{-1}(\varphi(\rho_1\cap\rho_2))\cap\rho_1$
is a simplex in ${\cal R}$.

Now consider two liftings $\tilde{x}$ and $g\cdot\tilde x$ of $x$ in $\U$, 
and their open star neighborhoods $\Sigma_{\tilde x}$ and
$\Sigma_{g\cdot \tilde x}$. Suppose 
$\Sigma_{\tilde x}\cap \Sigma_{g\cdot \tilde x}\neq\emptyset$, and let $y$ 
be a point in this intersection. Thus $y\in\inter\tau$ where $\tilde x\in\tau$
and $g\cdot\tilde x\in\tau$. But $\varphi\colon\U\rightarrow U$ is one-one on 
$\tau$, as shown in Lemma \ref{lift}. So $\tilde x=g\cdot\tilde x$ and 
$\Sigma_{\tilde x}=\Sigma_{g\cdot \tilde x}$.
This proves the lemma.

\begin{gevolg}\label{goodcov}
Let ${\cal M}=(M,{\cal U})$ be any orbifold with atlas ${\cal U}$.
There exists an atlas ${\cal V}$ for ${\cal M}$ such that
\begin{mylist}
\item
${\cal V}$ refines ${\cal U}$;
\item
For every chart $(\tilde{V},H,\psi)$ in ${\cal V}$, both 
$\tilde V$ and $V$ are contractible;
\item
The intersection of finitely many charts in ${\cal V}$ is either 
empty or again a chart in ${\cal V}$. 
\end{mylist}
\end{gevolg}

\paragraph{Proof.} For a good triangulation ${\cal T}$ of $(M,{\cal U})$ 
as above, consider a simplex $\sigma$. The open star $V=\St(\sigma)$ is 
contained in $U$ for some chart $(\U,G,\varphi)$ in ${\cal U}$. 
Let $x$ be a point in the interior of $\sigma$, so that $V=\Sigma_x$. 
Choose a lifting $\tilde x$ of $x$ in $\U$. Then, by Lemma \ref{stars} 
and Remark \ref{rems}(\ref{contract}), the map 
$\varphi\colon\Sigma_{\tilde x}\rightarrow\Sigma_x$ is part of a chart for 
${\cal M}$. The collection of all these charts for open stars 
$V=\St(\sigma)$ is the required atlas ${\cal V}$.

\paragraph{}Following the terminology of \cite[page 42]{BT}, 
we call a cover ${\cal V}$ of ${M}$ by charts as in Proposition~\ref{goodcov}
a {\it good cover} of the orbifold ${\cal M}$.

\subsection{Sheaf cohomology.}\label{coho}
Let ${\cal M}=(M,{\cal U})$ be an orbifold. Recall (\cite{MP}) that a 
sheaf on ${\cal M}$ is given by the following data:
\begin{mylist}
\item
For each chart $(\U,G,\varphi)$ in ${\cal U}$ an (ordinary) sheaf of 
abelian groups ${\cal A}_{\U}$ on $\U$;
\item
For each embedding $$\lambda\colon(\U,G,\varphi)\rightarrow(\tilde{V},H,\psi)$$
an isomorphism 
$${\cal A}(\lambda)\colon{\cal A}_{\U}\stackrel{\sim}{\longrightarrow}%
\lambda^*({\cal A}_{\tilde{V}}).$$
These isomorphisms are required to be functorial in $\lambda$; i.e.\ if $$\mu\colon(\tilde{V},H,\psi)\rightarrow(\tilde{W},K,\chi)$$ is 
another embedding, then the following square commutes:
$$
\diagram
{\cal A}_{\U}\rto^{{\cal A}(\lambda)}\dto_{{\cal A}(\mu\lambda)}& %
\lambda^*({\cal A}_{\tilde{V}})\dto^{\lambda^*{\cal A}(\mu)}\\
(\mu\lambda)^*({\cal A}_{\tilde{W}})\rdouble&\lambda^*\mu^*({\cal A}_{\tilde{W}})
\enddiagram
$$
where `=' denotes the canonical isomorphism; 
\item\label{action}
It follows that each ${\cal A}_{\U}$ is a $G$-equivariant sheaf on $\U$, 
see \cite{MP}.
\end{mylist}

\paragraph{} With the obvious notion of morphisms between sheaves, these 
sheaves form an abelian category $\mbox{\it Ab}({\cal M})$ with enough injectives.

Referring to the examples in Section \ref{one}, we remark that if the 
orbifold is defined from the action of a compact Lie group $L$ on a 
manifold $N$, then this category is (equivalent to) the category of 
$L$-equivariant sheaves on $N$. And if the orbifold is defined from a 
suitable foliation ${\cal F}$ on a manifold $N$, it is (equivalent to) 
the category of holonomy-invariant sheaves on $N$.

For a sheaf ${\cal A}$ as above, a {\it global section} $s$ of ${\cal A}$ 
is by definition a system of sections $s_{\U}\in\Gamma(\U,{\cal A}_{\U})$, 
one for each chart $\U$, and compatible in the sense that for each embedding $\lambda\colon\U\rightarrow\tilde{V}$, the identity 
${\cal A}(\lambda)(s_{\U})=\lambda^*(s_{\tilde{V}})$ holds. The group of all 
these global sections is denoted $\Gamma({\cal M},{\cal A})$.
This defines a functor $\Gamma\colon{\it Ab}({\cal M})\rightarrow{\it Ab}$, 
into the category of abelian groups, which is right exact and preserves injectives.
For an abelian sheaf ${\cal A}$, one defines
$$H^n({\cal M},{\cal A})=(R^n\Gamma)({\cal A}).$$
This definition of the sheaf cohomology of ${\cal M}$ is just a special case 
of the cohomology of a topos \cite{SGA42}, and hence it satisfies all the standard
functoriality and invariance properties. By way of example, we mention some
of these properties:

\subsubsection{Standard properties.}
\begin{mylist}
\item
Any (strong, \cite{MP}) map between orbifolds $f\colon{\cal N}\rightarrow{\cal M}$ 
induces an exact functor $f^*\colon{\it Ab}({\cal M})\rightarrow{\it Ab}({\cal N})$, 
and an induced homomorphism 
$f^*\colon H^n({\cal M},{\cal A})\rightarrow H^n({\cal N},f^*{\cal A})$.
\item
There is also an (adjoint) functor 
$f_*\colon {\it Ab}({\cal N})\rightarrow{\it Ab}({\cal M})$ and a corresponding 
(`Leray') spectral sequence
$$
E^{p,q}_{2}=H^p({\cal M},R^qf_*({\cal B}))\Rightarrow H^{p+q}({\cal N},{\cal B})
$$
see \cite[Expos\'e V, Section 5]{SGA42}.
\item
There are adjoint functors
$$
\pi_*\colon{\it Ab}({\cal M})\leftrightarrows{\it Ab}(M)\,\colon\negthinspace\pi^*,
$$
where ${\it Ab}(M)$ is the category of abelian sheaves on the underlying space $M$. 
There is a corresponding Leray spectral sequence
$$
E^{p,q}_2=H^p(M,R^q\pi_*({\cal A}))\Rightarrow H^{p+q}({\cal M},{\cal A}),
$$
where for a point $x\in M$, the stalk $R^q\pi_*({\cal A})_x$ is 
$H^q(G_{\tilde{x}},{\cal A}_{\tilde{x}})$. (Here ${\cal A}_{\tilde{x}}$
is the stalk of ${\cal A}_{\U}$ at any lifting $\tilde{x}\in\U$ of $x$ for some 
chart $(\U,G,\varphi)$, and the stabilizer $G_{\tilde{x}}$ acts on 
${\cal A}_{\tilde{x}}$ since ${\cal A}_{\U}$ is $G$-equivariant.
\item
(Mayer-Vietoris) If $M=U\cup V$ is a union of two open sets, the orbifold structure 
${\cal M}$ restricts to orbifold structures ${\cal M}|U$ and ${\cal M}|V$, and 
there is a Mayer-Vietoris sequence
$$
\rightarrow H^n({\cal M}|U,{\cal A}|U)\oplus H^n({\cal M}|V,{\cal A}|V)\rightarrow
H^n({\cal M}|U\cap V,{\cal A}|U\cap V)\rightarrow H^{n+1}({\cal M},{\cal A})\rightarrow 
$$
\item \label{cech-spec}
Define a presheaf ${\cal H}^q({\cal A})$ on the underlying space $M$, by 
${\cal H}^q({\cal A})(U)=H^q({\cal M}|U,{\cal A}|U)$. For any open cover ${\cal V}$
of $M$, there is a spectral sequence
$$
E^{p,q}_2=\check{H}^p({\cal V},{\cal H}^q({\cal A}))\Rightarrow H^{p+q}({\cal M, A})
$$
with as $E^{p,q}_2$-term the \v{C}ech cohomology of this cover ${\cal V}$
with coefficients in this presheaf ${\cal H}^p({\cal A})$; see 
\cite[Expos\'e V, Section 3]{SGA42}. 
\item
We single out a special case of this last property (\ref{cech-spec}). A sheaf on 
${\cal M}$ is said to be {\it locally constant} if each of the (ordinary) sheaves 
${\cal A}_{\U}$ is locally constant.
Now suppose ${\cal V}$ is a `good' atlas for ${\cal M}$, as in Corollary \ref{goodcov}, 
and let ${\cal A}$ be any locally constant sheaf.
For any chart $(\tilde{V},H,\psi)$ in ${\cal V}$, the restriction 
${\cal A}_{\tilde{V}}$ is a constant sheaf on a contractible space ${\tilde{V}}$,
and $H^q({\cal M}|V,{\cal A}|V)$ is the cohomology of the group $G=G_V$ with 
coefficients in the $G_V$-module $A_V=\Gamma(\tilde{V},{\cal A}_{\tilde{V}})$.
Thus, for a locally constant sheaf ${\cal A}$ and a good cover ${\cal V}$, the 
spectral sequence takes the form
$$
\check{H}^p({\cal V},V\mapsto{ H}^q(G_V,A_V))\Rightarrow H^{p+q}({\cal M,A}).
$$
\end{mylist}

\section{Simplicial complexes for orbifolds.}\label{simpcom}

Our purpose in this section is to describe explicitly for any orbifold 
${\cal M}=(M,{\cal U})$ a simplicial set $S$, with the property that any 
locally constant sheaf ${\cal A}$ on ${\cal M}$ induces a local system of 
coefficients $A$ on $S$ for which there is a natural isomorphism
$$
H^*({\cal M, A})\cong H^*(S,A),
$$
see Theorem~\ref{main} below.

\subsection{The simplicial set.}\label{simpset} Let us fix an $n$-dimensional 
orbifold ${\cal M}$ with underlying space $M$ and (chosen) atlas  ${\cal U}$; 
let us also fix a triangulation ${\cal T}$ of
${\cal M}$, and write $S_0$ for the set of $n$-simplices.
We assume that the triangulation is adapted to ${\cal U}$, as described in 
Section \ref{trian}. Recall that this means that ${\cal T}$ has the following 
properties:\\

\begin{mylist}
\item\label{i}
For each $n$-simplex $\sigma\in S_0$ there is a chart 
$(\U_\sigma,G_\sigma,\varphi_\sigma)$, such that $\sigma\subseteq U_\sigma$;
\item\label{ii}
For each simplex $\tau$ there is a face $\tau'\subseteq\tau$ such that 
the isotropy is constant on $\tau-\tau'$; in particular every simplex has 
a vertex
$v(\tau)$ with maximal isotropy.\\
\end{mylist}

\noindent We assume that a choice of charts $U_\sigma$ and vertices $v(\tau)$ 
as in (\ref{i}) and (\ref{ii}) above has been made.
(Note that we do not require the stronger property of being `good' for the 
triangulation, because in some examples that would force us to construct 
a simplicial complex which is bigger than necessary, and hence less 
suitable for calculations.)

We now construct the simplicial set $S=S({\cal T})$ with the same cohomology 
as ${\cal M}$. The description of $S$ will use various choices, besides the 
charts $U_\sigma$ and the vertices $v(\tau)$ already mentioned.
First of all, choose for each simplex $\sigma\in S_0$ a lifting $\tilde{\sigma}$ 
as in Lemma \ref{lift}, mapped homeomorphically to $\sigma\subseteq U_\sigma$ by $\varphi_\sigma\colon\U_\sigma\rightarrow U_\sigma$. Next, fix for each vertex 
$v$ of ${\cal T}$ a neighborhood $U_v$ of $v$ and a chart $(\U_v,G_v,\varphi_v)$ 
over $U_v$, so small that $U_v\subseteq U_\sigma$ whenever $v\in\sigma\in S_0$. 
Also fix a lifting $\tilde{v}\in\varphi_v^{-1}(v)$, and an embedding $\lambda_{\sigma,v}\colon\U_v\mono\U_\sigma$ with 
$\lambda_{\sigma,v}(\tilde{v})\in\tilde{\sigma}$.

We will not require these $U_v$ to belong to the original atlas ${\cal U}$.
In fact, they can be chosen so small that $G_v$ is also the isotropy 
group of $v$ in the chart $(\U_v,G_v,\varphi_v)$, so that the notation is 
unambiguous;
in this case the lifting $\tilde{v}$ is unique.

Let $\sigma_0,\sigma_1\in S_0$ be two simplices of maximal dimension 
$n$, and assume $\sigma_0\cap\sigma_1\neq\emptyset$. Below, in 
Section~\ref{five},
we will construct for any two simplices $\tau$ and $\rho$ with
\begin{equation}\label{subsimp}
\tau\subseteq\rho\subseteq\sigma_0\cap\sigma_1
\end{equation} 
an injective map
\begin{equation}\label{mu}
\mu_{\tau,\rho,\sigma_0,\sigma_1}\colon G_{v(\tau)}\rightarrow G_{v(\rho)}.
\end{equation}
We will write $\mu$ or $\mu_{\tau,\rho}$ if the (other) subscripts are clear
from the context. This map $\mu$ in (\ref{mu}) will not be a homomorphism 
in general; it will map $G_{v(\tau)}$ to a coset of a conjugate of the 
subgroup $G_{v(\tau)}=\{g\in G_{v(\rho)}|\,g\cdot v(\tau)=v(\tau)\}$ in 
$G_{v(\rho)}$. However the construction will have the following 
multiplicative property:
if 
$$
\tau\subseteq\rho\subseteq\sigma_0\cap\sigma_1\cap\sigma_2
$$
then, for $h_1,h_2\in G_{v(\tau)}$,
\begin{equation}\label{mult}
\mu_{\sigma_0,\sigma_1}(h_1)\cdot\mu_{\sigma_1,\sigma_2}(h_2)=%
\mu_{\sigma_0,\sigma_2}(h_1h_2).
\end{equation}
Moreover, if $\tau=\rho$ and $\sigma_0=\sigma_1$ then 
$\mu\colon G_{v(\tau)}\rightarrow G_{v(\tau)}$ is the identity.

With these choices made, the simplicial set $S=S({\cal T})$ can be described. 
As already defined above,
$$
S_0=\{\sigma|\,\sigma \mbox{ is an $n$-simplex of }{\cal T}\}.
$$
Furthermore, for $k\geq 1$,
$$
S_k=\sum_{\stackrel{\scriptstyle\sigma_0,\cdots,\sigma_k\in S_0}{ \sigma_0\cap\cdots\cap\sigma_k\neq\emptyset}} G^k_{v(\sigma_0\cap\cdots\cap\sigma_k)}.
$$
An element of $S_k$ can also be denoted by
\begin{equation}\label{nerve}
\sigma_0\stackrel{g_1}{\leftarrow}\sigma_1\leftarrow\cdots%
\stackrel{g_k}{\leftarrow}\sigma_k,
\end{equation}
to suggest the analogy with nerves. So $\sigma_0,\cdots,\sigma_k$ in 
(\ref{nerve}) are $n$-simplices, and $g_1,\cdots,g_k\in G_{v(\tau)}$ 
where $v(\tau)$ is the chosen vertex with maximal isotropy on 
$\tau=\sigma_0\cap\cdots\cap\sigma_k$. The degeneracy maps 
$s_i\colon S_{k-1}\rightarrow S_k$ $(i=0,\cdots,k-1)$ are defined in the usual way,
$$
s_i(\sigma_0\stackrel{g_1}{\leftarrow}\cdots\stackrel{g_{k-1}}{\leftarrow}%
\sigma_{k-1})=(\sigma_0\stackrel{g_1}{\leftarrow}\cdots\sigma_i\stackrel{1}%
{\leftarrow}\sigma_i\cdots\stackrel{g_{k-1}}{\leftarrow}\sigma_{k-1}).
$$
The face maps $d_j\colon S_k\rightarrow S_{k-1}$ $(j=0,\cdots,k)$ are 
defined by means of the maps $\mu$ in (\ref{mu}), as
$$
d_j(\sigma_0\stackrel{g_1}{\longleftarrow}\cdots\stackrel{g_k}%
{\longleftarrow}\sigma_k)=\negthinspace
\left\{\begin{array}{ll}
\negthinspace\negthinspace\negthinspace\sigma_1\stackrel{\mu(g_2)}%
{\longleftarrow}\sigma_2\longleftarrow\cdots\stackrel{\mu(g_k)}%
{\longleftarrow}\sigma_k, & \negthinspace\negthinspace\negthinspace(j=0)\\
\negthinspace\negthinspace\negthinspace\sigma_0\stackrel{\mu(g_1)}%
{\longleftarrow}\cdots\sigma_{j-1}\stackrel{\mu(g_jg_{j+1})}%
{\longleftarrow}\sigma_{j+1}\cdots\longleftarrow\sigma_k,%
&\negthinspace\negthinspace\negthinspace(0<j<k)\\
\negthinspace\negthinspace\negthinspace\sigma_0\stackrel{\mu(g_1)}%
{\longleftarrow}\cdots\stackrel{\mu(g_{k-1})}{\longleftarrow}%
\sigma_{k-1},&\negthinspace\negthinspace\negthinspace(j=k),
\end{array}\right.
$$
where the $\mu$'s carry the following subscripts.
Write $\tau=\sigma_0\cap\cdots\cap\sigma_k$ and $\rho=\sigma_0\cap\cdots\hat{\sigma}_j\cap\cdots\cap\sigma_k$.
Then
$$
\mu(g_j)=\mu_{\tau,\rho,\sigma_{j-1},\sigma_j}(g_j)
$$
and 
\begin{eqnarray*}
\mu(g_jg_{j+1})&=&\mu_{\tau,\rho,\sigma_{j-1},\sigma_{j+1}}(g_jg_{j+1})\\
&=&\mu_{\tau,\rho,\sigma_{j-1},\sigma_j}(g_j)\mu_{\tau,\rho,\sigma_{j},%
\sigma_{j+1}}(g_{j+1}),
\end{eqnarray*}
the latter by (\ref{mult}). The simplicial identities now follow easily.

Before we define these maps $\mu$, we state the theorem. Let ${\cal A}$ 
be a locally constant sheaf on ${\cal M}$. This sheaf induces in a natural 
way a local
system of coefficients $A$ on the simplicial set $S_\bullet$, with, 
for any $\sigma\in S_0$,
$$
A_\sigma=\Gamma(\tilde{\sigma},{\cal A}_{\tilde{U}_\sigma}).
$$
To describe the twisting, observe that, since ${\cal A}$ is locally 
constant, ${\cal A}_{\U_\sigma}$ is constant on $\tilde{\sigma}$. 
So for any vertex $\tilde{w}\in\tilde\sigma$ there is a canonical 
isomorphism from the stalk at $\tilde{w}$,
$$
\diagram
{\cal A}_{\U_\sigma,\tilde{w}}\rto^\sim&A_\sigma.
\enddiagram
$$
Modulo these isomorphisms, the twisting by an element 
$(\sigma_0\stackrel{g}{\leftarrow}\sigma_1)\in S_1$ is now defined as the 
dashed map in the diagram of maps between stalks
$$
\diagram
{\cal A}_{\U_{\sigma_1},\tilde v_1}\ddouble\xdashed[0,3]|>\tip&&&%
{\cal A}_{\U_{\sigma_0},\tilde v_0}\ddouble\\
\lambda_1^*({\cal A}_{\U_{\sigma_1}})_{\raisebox{-1.5pt}%
{${\scriptstyle {\tilde{v}}}$}}&{\cal A}_{\U,v}\lto\rto^{\overline{g}}&%
{\cal A}_{\U,v}\rto&\lambda_0^*({\cal A}_{\U_{\sigma_0}})_{\raisebox{-1.5pt}%
{${\scriptstyle {\tilde{v}}}$}}
\enddiagram
$$
Here $v=v(\sigma_0\cap\sigma_1)$, so $g\in G_v$. Furthermore, $\lambda_i=\lambda_{\sigma_i,v}\colon\U_v\rightarrow\U_{\sigma_i}$ is the 
chosen 
embedding $(i=0,1)$, with the property that 
$\tilde{v}_i=\lambda(\tilde{v})\in\tilde{\sigma}_i$.
Finally, $\overline{g}$ in the diagram denotes the left action by 
$g\in G_{\tilde{v}}$ (cf.\ condition (\ref{action}) in the description 
of sheaves in Section \ref{coho}).
 
\begin{stel}\label{main}
For any triangulated orbifold ${\cal M}$ with associated simplicial set
$S$ as above, and for any locally constant sheaf ${\cal A}$ on ${\cal M}$,
there is a natural isomorphism
$$
H^{\raisebox{.8pt}{${\scriptscriptstyle \,\bullet}$}}({\cal M,A})\cong H^{\raisebox{.8pt}{${\scriptscriptstyle \,\bullet}$}}(S,A).
$$
\end{stel}

We will now first define the maps $\mu$ involved in the definition of $S$. 
The proof of the theorem will be given in Section~\ref{proof}. 

\subsection{Construction of $\mu$.}\label{five}
Fix $\sigma_0,\sigma_1\in S$ and
$\tau\subseteq\rho\subseteq\sigma_0\cap\sigma_1$ as in (\ref{subsimp}). 
Write $v=v(\tau)$ and $w=v(\rho)$ for the corresponding vertices with 
maximal isotropy on $\tau$ and $\rho$, respectively.
We will construct
\begin{equation}\label{muvw}
\mu\colon G_v\rightarrow G_w.
\end{equation}
Let $\theta$ be the 1-simplex joining $v$ and $w$ in $\rho$. Then the 
isotropy
group of any interior point of $\theta$ agrees with that of $v$, while 
that of $w$ is possibly larger, $w$ being maximal on $\theta$. Recall 
that we have already chosen
$$
\begin{array}{ll}
\lambda_i=\lambda_{\sigma_i,v}\colon\U_v\rightarrow\U_{\sigma_i},&%
\lambda_i(\tilde{v})\in\tilde\sigma_i,\\
\chi_i=\lambda_{\sigma_i,w}\colon\U_w\rightarrow\U_{\sigma_i},&%
\chi_i(\tilde{w})\in\tilde\sigma_i.
\end{array}
$$
First consider the special case that there is a chart 
$(\tilde W,H,\psi)$ for a neighborhood $W\supseteq \theta$, for which 
$U_v,U_w\subseteq W$ and $W\subseteq U_{\sigma_0}\cap U_{\sigma_1}$. 
(Such a chart exists, for example, 
when one starts with a `good' atlas ${\cal U}$ as in Section \ref{trian}, 
in which case one can take $W=U_{\sigma_0}\cap U_{\sigma_1}$.) 
Choose an embedding
$$
\gamma_0\colon\tilde{W}\rightarrow\U_{\sigma_0}
$$
with
$\gamma_0(\tilde{W})\cap\tilde\sigma_0\neq\emptyset$.
Let $\tilde\theta=\gamma_0^{-1}(\tilde\sigma_0\cap\varphi_{\sigma_0}^{-1}%
(\theta))\subseteq\tilde W$ be a lifting of $\theta$, and choose the 
other embedding
$$
\gamma_1\colon\tilde{W}\rightarrow\U_{\sigma_1}
$$
such that $\gamma_1(\tilde\theta)\subseteq\tilde\sigma_1$.
Also choose embeddings
$$
\tilde{U}_v\stackrel{\alpha}{\longrightarrow}\tilde{W}\stackrel{\beta}%
{\longleftarrow}\U_w
$$
with $\alpha(\tilde{v}), \,\beta(\tilde{w})\in\tilde\theta$. 

Now observe that, since $\gamma_0\alpha$ and $\lambda_0$ both map 
$\tilde{v}$
into $\tilde\sigma_0$, one has $\gamma_0\alpha(\tilde v)=\lambda_0(\tilde v)$, 
and hence there is a $g_0\in G_v$ such that $\gamma_0\alpha=\lambda_0 g_0$ 
(cf.\ Section \ref{one}, Remark (\ref{two})).
Similarly, we find $g_1\in G_v$ and $h_0,h_1\in G_w$, such that 
$$
\gamma_i\alpha=\lambda_i g_i,\quad\gamma_i\beta =\chi_i h_i,\quad(i=0,1).
$$

We claim that for any $k\in G_v$ there is a (unique) $m\in G_w$ such that
\begin{equation}\label{m}
\alpha(g_0^{-1}kg_1)=\beta(h_0^{-1}mh_1).
\end{equation}
To see this, recall first that $H$ denotes the group of the chart $\tilde{W}$,
and consider
$$
\ell =\alpha(g_0^{-1}kg_1).
$$
Note that $\ell $ fixes $\alpha(\tilde{v})\in\tilde\theta$.
Since the isotropy along $\alpha(\tilde\theta)$ does not decrease from
$\alpha(\tilde v)$, it follows that $\ell $ fixes $\tilde\theta$, and hence 
also the point $\beta(\tilde w)$. Thus $\ell =\beta(\ell')$ for some 
$\ell'\in G_w$.
Now let
$$
m=h_0 \ell'h_1^{-1}.
$$
Then $\alpha(g_0^{-1}kg_1)=\beta(h_0^{-1}mh_1)$, as required for (\ref{m}).

We now define $\mu$ in (\ref{muvw}) by
\begin{equation}\label{mukl}
\mu(k)=m.
\end{equation}
It can be shown that this definition is independent of the choices made above 
(of $\tilde W$, $\alpha$, $\beta$).

\begin{figure}[htb]
$$
\makebox[0pt]{\epsfxsize=.8\textwidth\epsffile{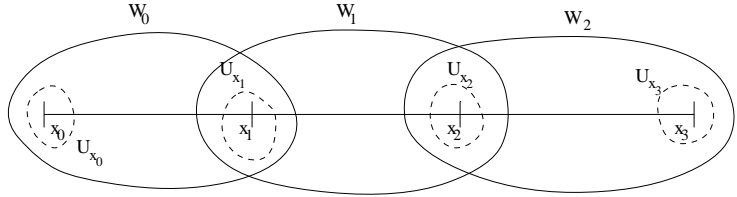} } 
$$
\caption{\label{chain} A subdivision of $\theta$ with $n=3$.}
\end{figure}

\paragraph{}This defines $\mu$ in the special case, where the 1-simplex
$\theta$ is contained in a chart $W$ with $U_{\sigma_0}\cap U_{\sigma_1}%
\supseteq W\supseteq U_v,U_w$. In the general case, choose a subdivision 
of $\theta$ into smaller 1-simplices, with vertices
$$
v=x_0,\,x_1,\cdots,\, x_n=w.
$$
Choose charts
$$
\U_v=\U_{x_0},\,\U_{x_1},\cdots,\,\U_{x_n}=\U_w,
$$
liftings $\tilde{x}_i\in\U_{x_i}$, and embeddings
$$
\lambda_{x_j,i}\colon\U_{x_j}\rightarrow\U_{\sigma_i}\quad(i=0,1;\,j=0,\cdots,n),
$$
with $\lambda_{x_j,i}(\tilde{x}_j)\in\tilde{\sigma}_i$, and coinciding with 
the embeddings already chosen for $j=0,n$. Furthermore, choose this subdivision
of $\theta$ sufficiently fine, and these charts $\U_{x_i}$ sufficiently small,
so that each 1-simplex between $x_i$ and $x_{i+1}$ is contained in a chart 
$W_i\supseteq U_{x_i}, U_{x_{i+1}}$, as in Figure \ref{chain}.
Now define $\mu_k\colon G_{x_k}\rightarrow G_{x_{k+1}}$ exactly as the maps
$\mu\colon G_v\rightarrow G_w$ defined in (\ref{mukl}) above, for 
$k=0,\cdots,n-1$, and let 
$$
\mu\colon G_v=G_{x_0}\rightarrow G_{x_n}=G_w
$$
be the composition
$$
\mu=\mu_{n-1}\circ\cdots\circ\mu_1\circ\mu_0.
$$
It can be shown that this definition of $\mu$ is again independent of the 
various choices. In particular, for a finer subdivision of $\theta$ and a 
refinement
of the system of open sets $W_j$ and $U_{x_i}$, one obtains the same map $\mu$. 
We omit the details.

\section{Example: the teardrop orbifold.}\label{examples}

In this section we will apply the construction of the simplicial set 
$S$ from the previous 
sections to the teardrop orbifold (as described in \cite{Thurston}), and 
calculate its cohomology groups.
As before ${\cal A}$ denotes a locally constant sheaf of coefficients on the
orbifold considered, and $A$ denotes the induced local system of coefficients
on $S$.

\subsection{The triangulation.} 
\begin{figure}[htb]
$$
\makebox[0pt]{\epsfysize=.12\textheight\epsffile{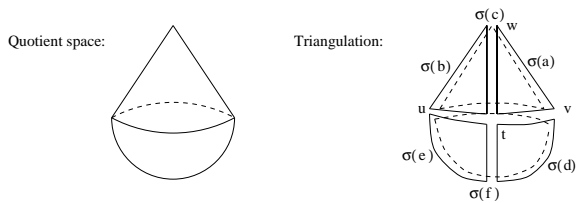}}
$$
\caption{\label{teardrop} The quotient space of the teardrop orbifold 
with a triangulation.}
\end{figure}

The quotient space of the teardrop orbifold is the 2-sphere with one cone
point of order $n$. A chart around the cone point consists of an 
open disk $\U$ in ${\mathbb R}^2$ with structure group $C_n$ (the finite 
cyclic group of order $n$), which acts on $\U$ by rotations.
We will denote this orbifold by {\it Tear-n}.
Figure \ref{teardrop} shows a picture of the quotient space and a  triangulation
of this orbifold.
The simplices $\sigma(a),\,\sigma(b),\,\sigma(d),$ and $\sigma(e)$
are on the front and the simplices $\sigma(c)$ and $\sigma(f)$ are on 
the back of the teardrop. Moreover, $t=\sigma(abde)$, $u=\sigma(bcef)$, 
$v=\sigma(acdf)$ and  $w=\sigma(abc)$, where 
$\sigma(i_0\cdots i_n)=\sigma(i_0)\cap\cdots\cap\sigma(i_n)$.

\begin{figure}[htb]
$$
\makebox[0pt]{\epsfysize=.14\textheight\epsffile{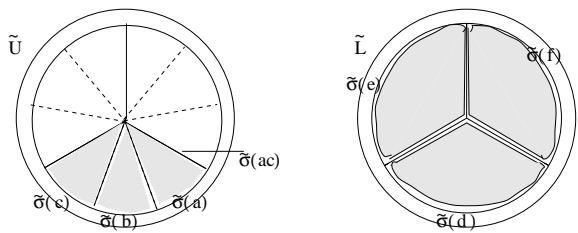}}
$$
\caption{\label{liftings} Chosen liftings for the simplices of 
Figure \ref{teardrop}, when $n=3$.}
\end{figure}

The atlas we use for this orbifold consists of eight charts: an open disk 
$\tilde L$ with a trivial structure group to cover the lower half of the 
quotient space, and an open disk $\U$ with structure group $C_n$ (acting 
by rotations) to cover the upper half, and six charts to cover the equator 
in order to satisfy the compatibility condition for atlases. Note that the 
triangulation is adapted to this atlas.
Figure \ref{liftings} shows the liftings (in the charts) of the 1- and 
2-simplices in the triangulation, as needed for the construction of the
simplicial representation for this orbifold, for the case that $n=3$.
The simplices $\sigma(a),\,\sigma(b),$ and $\sigma(c)$ are subsets of
$U$; their liftings, denoted by $\tilde\sigma(-)$,  are shown in the 
left hand chart $\U$. Similarly the liftings of $\sigma(d),\,\sigma(e)$, 
and $\sigma(f)$
are shown in the right-hand chart $\tilde L$.  
So $U_{\sigma(a)}=U_{\sigma(b)}=U_{\sigma(c)}=\U$, whereas 
$U_{\sigma(d)}=U_{\sigma(e)}=U_{\sigma(f)}=\tilde L$.

\subsection{The $\mu$-maps.} Using these liftings we define the maps
$$
\mu_{\sigma(ijk),\sigma(ij),\sigma(i),\sigma(j)}\colon %
G_{v(\sigma(ijk))}\rightarrow G_{v(\sigma(ij))}
$$
for every pair $\sigma(i)$ and $\sigma(j)$ for which the intersection 
is a 1-simplex $\sigma(ij)$ and $\sigma(ijk)$ is a vertex with the same 
isotropy
group as the interior of $\sigma(ij)$.
Note that $G_{v(\sigma(ijk))}$ (where $\sigma(i)$, $\sigma(j)$ and 
$\sigma(k)$
are three different simplices) is trivial except for 
$G_w=G_{v(\sigma(abc))}=C_n$.
So the only possibly non-trivial $\mu$-maps are those with codomain 
$G_{w}=G_{v(\sigma(ab))}=G_{v(\sigma(bc))}=G_{v(\sigma(ac))}$.
In order to construct
$$
\mu_{t,\sigma(ab),\sigma(a),\sigma(b)}\colon G_{t}\rightarrow G_{v(\sigma(ab))}=G_{w},
$$
note that $\theta$, as used in Section~\ref{five}, is $\sigma(a)\cap\sigma(b)$.
We choose $W=U$ and $\tilde\theta=\tilde\sigma(a)\cap\tilde\sigma(b)$, 
so both $\gamma_0$ and $\gamma_1$, as used in Section~\ref{five}, 
are the identity map. Moreover, we choose 
\begin{equation}\label{een}
\lambda_{\sigma(a),w}=\lambda_{\sigma(b),w}
\end{equation}
(notation as in Section~\ref{simpset}),
whereas $\lambda_{\sigma(a),t}$ and $\lambda_{\sigma(b),t}$ have 
to be the same, because the liftings
of $t$ in $\tilde{\sigma}(a)$ and $\tilde{\sigma}(b)$ are the same point in $\U$. 
And it is obvious that
$$
\mu_{t,\sigma(ab),\sigma(a),\sigma(b)}(1_{G_{t}})=1_{G_{v(\sigma(ab))}}.
$$
Also (by the multiplicative property (\ref{mult}) of $\mu$),
$$
\mu_{t,\sigma(ab),\sigma(b),\sigma(a)}(1_{G_{t}})=1_{G_{v(\sigma(ab))}}.
$$
For similar reasons, when we choose 
\begin{equation}\label{twee}
\lambda_{\sigma(b),w}=\lambda_{\sigma(c),w}
\end{equation}
we find that both maps 
$$\mu\colon G_{u}\rightarrow G_{v(\sigma(bc))}$$
are described by 
$$\mu(1_{G_{u}})=1_{G_{v(\sigma(bc))}}.$$
The only non-trivial $\mu$-maps are 
$\mu_{v,\sigma(ac),\sigma(a),\sigma(c)}$, and 
$\mu_{v,\sigma(ac),\sigma(c),\sigma(a)}$ from $G_{v}$
into $G_{v(\sigma(ac))}=G_{w}$. 
With the notation of Section \ref{simpcom}, $\theta=\sigma(ac)$ and
we choose $W=U_{v}=U_{w}=U$ with 
$\tilde\theta=\tilde\sigma(ac)\subset\tilde\sigma(a)$ as in Figure~\ref{liftings}. 
Moreover, we take 
$$
\lambda_{\sigma(a),w}=\lambda_{\sigma(a),v}=\gamma_1=\id.
$$
Then (\ref{een}) and (\ref{twee}) above induce that 
$\lambda_{\sigma(c),w}=\id.$ It follows from the choice of the liftings 
in Figure~\ref{liftings}, that 
$$
\gamma_0=\lambda_{\sigma(c),v}=\rho,
$$
where $\rho$ is rotation over $2(n-1)\pi/n$, which generates $G_{w}$. 
We conclude that the group elements $g_0, g_1,$ and $h_1$ are all the 
identity element, and $h_0=\rho$, so: 
$$
\mu_{v,\sigma(ac),\sigma(c),\sigma(a)}(1)=\rho.
$$
Similarly one can show that:
$$
\mu_{v,\sigma(ac),\sigma(a),\sigma(c)}(1)=\rho^{n-1}.
$$
(This also follows from the fact that they have to be each other's 
inverses by the multiplicative property of $\mu$ as stated in 
(\ref{mult}).)

\subsection{The simplicial complex.} With these $\mu$-maps the simplicial 
complex can be described as follows.
\begin{eqnarray*}
S_0&=&\{\sigma(a),\sigma(b),\sigma(c),\sigma(d),\sigma(e),\sigma(f)\},\\
S_n&=&\{(\sigma(i_0\cdots i_n),g_1,\cdots,g_n);\,g_1,\cdots,g_n\in%
 G_{\sigma(i_0\cdots i_n)}\mbox{ and }\sigma(i_0\cdots i_n)\neq\emptyset\}. 
\end{eqnarray*} 
The degeneracy maps $s_i$ are trivial as described in Section \ref{simpset},
and the face maps are straightforward compositions on all combinations
which do not contain neighboring $a$ and $c$ in their $\sigma$-part, since 
the $\mu$-maps are trivial in these cases.
For the $ac$-combinations, we have to use the non-trivial part
of the $\mu$-maps described above. For example, deleting the sigma's 
from the notation, we can describe the face operators on $S_2$ by: 
\begin{eqnarray*}
d_0(f\stackrel{1}{\longleftarrow} c\stackrel{1}{\longleftarrow}a)&=&%
(c\stackrel{\rho}{\longleftarrow}a),\\
d_1(f\stackrel{1}{\longleftarrow} c\stackrel{1}{\longleftarrow}a )%
&=&(f\stackrel{1}{\longleftarrow}a),\\ 
d_2(f\stackrel{1}{\longleftarrow} c\stackrel{1}{\longleftarrow}a)%
&=&(f\stackrel{1}{\longleftarrow}c).
\end{eqnarray*}
So for $\alpha\in A_{\sigma(a)}$,
we find that $(f\stackrel{1}{\longleftarrow}c) \cdot( %
(c\stackrel{\rho}{\longleftarrow}a)\cdot\alpha)= %
(f\stackrel{1}{\longleftarrow}a)\cdot\alpha$.
We use this in the following calculation, 
\begin{eqnarray*}
(a\stackrel{\,\rho}{\leftarrow}a)\cdot\alpha&=&(a\stackrel{1}{\leftarrow}c)%
\cdot((c\stackrel{\,\rho}{\leftarrow}a)\cdot\alpha)\\
&=&(a\stackrel{1}{\leftarrow}b)\cdot((b\stackrel{1}{\leftarrow}c)\cdot%
((c\stackrel{\,\rho}{\leftarrow}a)\cdot\alpha))\\
&=&(a\stackrel{1}{\leftarrow}d)\cdot((d\stackrel{1}{\leftarrow}b)\cdot%
((b\stackrel{1}{\leftarrow}f)\cdot((f\stackrel{1}{\leftarrow}c)\cdot%
((c\stackrel{\,\rho}{\leftarrow}a)\cdot\alpha))))\\
&=&(a\stackrel{1}{\leftarrow}f)\cdot((f\stackrel{1}{\leftarrow}d)\cdot%
((d\stackrel{1}{\leftarrow}e)\cdot((e\stackrel{1}{\leftarrow}b)\cdot%
((b\stackrel{1}{\leftarrow}f)\cdot\\
&&\cdot((f\stackrel{1}{\leftarrow}a)\cdot\alpha)))))\\
&=&(a\stackrel{1}{\leftarrow}f)\cdot((f\stackrel{1}{\leftarrow}e)%
\cdot((e\stackrel{1}{\leftarrow}f)\cdot((f\stackrel{1}{\leftarrow}a)%
\cdot\alpha)))\\
&=&(a\stackrel{1}{\leftarrow}f)\cdot((f\stackrel{1}{\leftarrow}a)\cdot\alpha)\\
&=&\alpha
\end{eqnarray*}
(This calculation illustrates the effect of the fact that a loop around 
the singular point $v(\sigma(abc))$ is contractible via the lower half of the teardrop.)

\subsection{The cohomology groups} We conclude from the above calculation 
that $C_n$ acts trivially on $A$, and it is not difficult to derive that:
$$
H^0(\mbox{\it Tear-n};{\cal A})=A\mbox{ and }H^1(\mbox{\it Tear-n};{\cal A})=0.
$$

To find the second cohomology group, we have to do some cocycle-coboundary 
calculations in the local system of coefficients $A$ on the simplicial set 
$S$. A cohomology class in $H^2$ is represented by a cocycle $\alpha\in Z^2(S,A)$, 
which we write as an $A$-valued map $\alpha\colon S_2\rightarrow A$. Since the 
second cohomology group of the plane is zero, we can choose the representant 
$\alpha$ in such a way that $\alpha(\sigma(ijk),g_1,g_2)=0$ when 
$\{i,j,k\}\cap\{d,e,f\}\neq\emptyset$ and moreover, \\
\begin{tabular}{rrr}
$\alpha(\sigma(aca),\rho^{n-1},\rho)=0$,& $\alpha(\sigma(cac),\rho,\rho^{n-1})=0$,&
$\alpha(\sigma(aba),1,1)=0$,\\
$\alpha(\sigma(bab),1,1)=0$,& $\alpha(\sigma(bcb),1,1)=0$,& $\alpha(\sigma(cbc),1,1)=0$.
\end{tabular}\\
It follows then from the cocycle-relations, that $\alpha$
is determined by the following data
\begin{equation}\label{data}
\left\{\begin{array}{lll}
\alpha(\sigma(caa),1,\rho^i)&\alpha(\sigma(bba),\rho^i,1)& \alpha(\sigma(cbb),1,\rho^i)\\
&\alpha(\sigma(aab),\rho^i,1)&\alpha(\sigma(ccb),\rho^i,1)\\
\alpha(\sigma(acc),1,\rho^i)&\alpha(\sigma(baa),1,\rho^i)& \alpha(\sigma(bbc),\rho^i,1)\\
\alpha(\sigma(aaa),\rho^i,\rho^j)&\alpha(\sigma(abc),1,1)&
\end{array}\right.
\mbox{ for }i,j\in\{1,2,\cdots,n-1\}.
\end{equation}
By choosing an appropriate coboundary, one can find an equivalent cocycle 
$\alpha'$, such that of all the values in (\ref{data}) above, only
$\alpha'(\sigma(abc),1,1)$ need not be zero.
Moreover, two cocycles $\alpha_1$ and $\alpha_2$ of this form are 
equivalent iff 
$\alpha_1(\sigma(abc),1,1)=\alpha_2(\sigma(abc),1,1)$.
We conclude:
$$
H^2(\mbox{\it Tear-n};{\cal A})=A.
$$

The higher degree cohomology groups of {\it Tear-n} can be calculated 
using 
the Mayer-Vietoris sequence for the upper and lower half-sphere, 
denoted by
$U_n$ and $L$, for $m>2$:
$$
{\spreaddiagramcolumns{-1.4pc}
\diagram
\rto & H^{m-1}(S^1;A)\ddouble\rto & H^m(\mbox{{\it Tear-n}};{\cal A})\rto\ddouble%
 & H^m(U_n;{\cal A})\oplus H^m(L;{\cal A})\rto\ddouble & H^m(S^1;{\cal A})\ddouble\rto & \\
\rto & 0\rto & H^m(\mbox{{\it Tear-n}};{\cal A})\rto & H^m(C_n;A)\oplus 0\rto & 0\rto & 
\enddiagram}
$$ 
We conclude:
\begin{stel}
The teardrop orbifold {\it Tear-n} has the following cohomology groups:
\begin{eqnarray*}
H^0(\mbox{{\it Tear-n}};{\cal A})&=&A;\\
H^1(\mbox{{\it Tear-n}};{\cal A})&=&0;\\
H^2(\mbox{{\it Tear-n}};{\cal A})&=&A;\\
H^m(\mbox{{\it Tear-n}};{\cal A})&=&H^m(C_n;A),\mbox{ for }m>2.
\end{eqnarray*}
\end{stel}

\section{Topological groupoids for orbifolds.}\label{seven}
Again, we fix an orbifold ${\cal M}$ with underlying space $M$ and 
atlas ${\cal U}$. Recall from Section \ref{coho} the category 
${\it Ab}({\cal M})$ of all abelian sheaves on ${\cal M}$. In 
\cite[Theorem 4.1]{MP} we proved that this category can be represented 
as the category of equivariant sheaves in various ways. Here we single 
out one particular such representation.

\subsection{Etale groupoids.} Let $G$ be a topological groupoid. 
As in loc. cit., we write $G_0$ for the space of objects and $G_1$ for 
the space of arrows, while the structure maps are denoted:
$$
\diagram
G_1\times_{G_0}G_1\rto^-{m} & G_1 \rto<.5ex>^{s}\rto<-.5ex>_{t} & G_0 %
\rto^u & G_1 \rto^i & G_1,
\enddiagram
$$
for composition, source, target, units and inverse, respectively.
As usual, we write $1_x$ for $u(x)$, $g^{-1}$ for $i(g)$, $g\circ h$ or 
$gh$ for $m(g,h)$, and $g\colon x\rightarrow y$ if $s(g)=x$ and $t(g)=y$. 
The groupoid $G$ is said to be {\it \'etale} if $s$ and $t$ are local 
homeomorphisms, and {\it proper} if $(s,t)\colon G_1\rightarrow G_0\times G_0$ 
is a proper map.

A $G$-{\it sheaf}  is a sheaf $A$ on the space $G_0$, equipped with a continuous
(say, right) action by $G_1$.
The category of all $G$-sheaves is denoted ${\it Ab}(G)$. We recall from \cite{MP}:

\begin{stel}\label{old}
For any orbifold ${\cal M}$ there exists a proper and \'etale topological 
groupoid $G$, for which there is a natural equivalence 
${\it Ab}({\cal M})\cong{\it Ab}(G)$.
\end{stel}

One construction of $G$ from ${\cal M}$, suggested in \cite{hae} and 
different from that in \cite{MP}, is the following. Let $G_0$ be the 
space of pairs
$(\tilde{x},\U)$ with $\tilde{x}\in\U\in{\cal U}$, topologized as the 
disjoint sum of the sets $\U$ in ${\cal U}$. An arrow $g\colon(\tilde{x},\U)\rightarrow(\tilde{y},\tilde{V})$ is an equivalence 
class of triples
$$
g=[\lambda,\tilde z, \mu]\colon\,\U\stackrel{\lambda}{\longleftarrow}\tilde W\stackrel{\mu}{\longrightarrow}\tilde{V},
$$
where $\tilde z\in\tilde W$ and $\lambda(\tilde z)=\tilde x$, 
$\mu(\tilde z)=\tilde y$. Here $\tilde W$ is another chart for ${\cal M}$, 
and $\lambda,\,\mu$ are embeddings. The equivalence relation is generated by
$$
[\lambda,\tilde{z},\mu]=[\lambda\nu,\tilde z',\mu\nu],
$$
for $\lambda,\,\tilde z,\,\mu$ above and $\nu\colon\tilde W'\rightarrow \tilde W$ 
another embedding, with $\nu(\tilde z')=\tilde z$. There is a natural topology 
on the set $G_1$ of all these equivalence classes, for which the source
and target maps $s,t\colon G_1\rightarrow G_0$ are each \'etale and together give
a proper map $G_1\rightarrow G_0\times G_0$. (See \cite{thesis} for details.)
Note that one can take any chart, not necessarily from the chosen atlas 
${\cal U}$, to represent an arrow, as it will always be equivalent to one 
represented by a chart in ${\cal U}$. (In fact, given a chart $\tilde W$ 
around $\varphi(\tilde x)$ which is embeddable into $\U$ and $\tilde V$, 
every arrow $g\colon(\tilde{x},\U)\rightarrow(\tilde{y},\tilde{V})$ can be 
represented by 
a triple involving $\tilde W$.)

For $G$ constructed in this way, it is not difficult to see that there is
an equivalence of categories ${\it Ab}({\cal M})\simeq{\it Ab}(G)$.

\paragraph{}
If $G$ is any topological groupoid, its nerve $\mbox{Nerve}(G)$ is the 
simplicial space whose $n$-simplices are strings $\stackrel{\rightarrow}{g}=(x_0\stackrel{g_1}{\leftarrow}\cdots%
\stackrel{g_n}{\leftarrow}x_n)$, 
equipped with the natural (fibered product) topology. One writes $G_n$ for 
the space $\mbox{Nerve}(G)_n$ of these $n$-simplices. If $A$ is any $G$-sheaf, 
it induces a sheaf $A^{(n)}$ on $G_n$, with stalk 
$A^{(n)}_{\stackrel{\rightarrow}{g}}=A_{x_n}$.

\begin{prop}\label{spec}
Let $G$ be any \'etale topological groupoid representing the orbifold 
${\cal M}$
as in Theorem \ref{old}. For each abelian sheaf ${\cal A}$ on ${\cal M}$, 
there is a natural spectral sequence
$$
E^{p,q}_2=H^pH^q(G_\bullet,A^{(\bullet)})\Rightarrow H^{p+q}({\cal M},{\cal A}).
$$
\end{prop}
(Here ${\cal A}$ corresponds by Theorem \ref{old} to a 
$G$-sheaf $A$, with induced sheaf $A^{(p)}$ on $G_p$, so that for fixed 
$q$, $H^q(G_\bullet,A^{(\bullet)})$ is a cosimplicial group.)

\paragraph{Proof.}
Using the equivalence of Theorem \ref{old}, this spectral sequence is simply 
a special case of the standard one,
$$
H^p H^q (G_{\bullet},A^{\bullet})\Rightarrow H^{p+q}(G,A),
$$
for \'etale topological groupoids (see \cite[V, (7.4.0.3)]{SGA42}). For later
use, we recall that the latter spectral sequence is constructed from the 
double complex $\Gamma(G_p,I^{q(p)})$, where 
$A\rightarrow I^0\rightarrow I^1\rightarrow\cdots$ is any injective resolution in 
${\it Ab}(G)$. This induces a resolution 
$A^{(p)}\rightarrow I^{0(p)}\rightarrow I^{1(p)}\rightarrow %
I^{2(p)}\rightarrow\cdots$ 
of sheaves on $G_p$ which is again injective.

\subsection{Proper groupoids.} The groupoid $G$ for ${\cal M}$ in 
Theorem \ref{old} is not unique.
(However, it is unique up to weak, or `Morita' equivalence \cite{MP}.)
For the proof of Theorem \ref{main}, we will use the following construction.
Let $H_0\mono G_0$ be a closed subspace. Then $H_0$ is the space of objects
of a topological groupoid $H$, with $H_1$ constructed as the fibered product
\begin{equation}\label{fibprod}
\diagram
H_1\rto|<\hole|<<\ahook\dto & G_1\dto\\
H_0\times H_0 \rto|<\hole|<<\ahook & G_0\times G_0.
\enddiagram
\end{equation}
In other words, $H$ is the full subgroupoid of $G$ on the space of 
objects $H_0\subseteq G_0$. (This topological groupoid $H$ is in general 
not \'etale.)

\begin{lemma}\label{subgr}
Assume the closed subspace $H_0\subseteq G_0$ has the property that the map:
$$
s\circ\pi_2\colon H_0\times_{G_0}G_1\rightarrow G_0,\quad(x\in H_0,%
(g\colon y\rightarrow x)\in G_1)\mapsto y
$$
is a proper surjection. Then the inclusion of groupoids $H\subseteq G$ 
induces an equivalence of categories ${\it Ab}(H)\simeq{\it Ab}(G)$.
\end{lemma}
 
\paragraph{Proof.} Standard. (See \cite{thesis} for details.)

\begin{gevolg}\label{subspec}
Let ${\cal M}$ be an orbifold, let $G$ be an \'etale groupoid representing
${\cal M}$ as in Theorem \ref{old}, and let $H\subseteq G$ be any subgroupoid
as constructed above. Then there exists a natural spectral sequence
$$
E^{p,q}_2=H^pH^q(H_\bullet, A^{(\bullet)}|H_\bullet)\Rightarrow H^{p+q}({\cal M,A}).
$$
\end{gevolg}
(Here $H_p\subseteq G_p=\mbox{Nerve}(G)_p$ and $A^{(p)}|H_p$ is the restricted sheaf.)

\paragraph{Proof.}
This is proved in the same way as in Proposition \ref{spec}, using that the 
injective sheaves $I^{(q)p}$ on $G_p$ restrict to soft sheaves on the closed 
subspace $H_p\subseteq G_p$.

\subsection{Proof of Theorem \ref{main}.}\label{proof}
Our purpose in this section is to construct a specific subgroupoid 
$H\subset G$ for ${\cal M}$
from a given triangulation. First, observe the following method for 
constructing a groupoid $H\subseteq G$ as in Lemma~\ref{subgr}. Let 
${\cal F}=\{F_i\}_{i\in I}$ be a locally finite cover of $M$ by compact sets, 
which refines the cover of charts, say $F_i\subset U_i$. Suppose there are 
chosen liftings $\tilde{F}_i\subset\U_i$, where the quotient map 
$\varphi_i\colon\U_i\rightarrow U_i$ of the chart maps $\tilde{F}_i$ 
homeomorphically to $F_i$. Let $H_0=\{(x,\tilde{F}_i)|x\in\tilde{F}_i,i\in I\}$ 
be the disjoint sum of these sets $\tilde{F}_i$. Then $H_0$ is a closed subset
of $G_0$. The induced groupoid $H$, with $H_1$ constructed as the fibered 
product in (\ref{fibprod}), is denoted $H({\cal F})$.

\begin{lemma}\label{subeq}
For any locally finite cover ${\cal F}$ as above, the map 
$s\pi_2\colon H_0\times_{G_0}G_1\rightarrow G_0$ is a proper 
surjection, and hence
the inclusion $H({\cal F})\subseteq G$ induces an equivalence 
${\it Ab}(H)\simeq{\it Ab}(G)\simeq {\it Ab}({\cal M})$. 
\end{lemma}

\paragraph{Proof.}
The second statement follows from the first by Lemma \ref{subgr}. 
It follows immediately from the fact that the $F_i$ cover $M$ that 
$s\pi_2\colon H_0\times_{G_0}G_1\rightarrow G_0$ is a surjection. 
To check that $s\pi_2$
is also  proper, it suffices to prove that each point in $G_0$ has a 
neighborhood $\tilde{V}$ over which $s\pi_2$ is a proper map 
$s\pi_2^{-1}(\tilde{V})\rightarrow\tilde{V}$. Choose $(\tilde{x},\U_j)\in G_0$
and let $\tilde{V}\subseteq\tilde{U}_j$ be a neighborhood of $\tilde{x}$ 
such that $V=\varphi_j(\tilde{V})\subseteq M$ meets only finitely many $F_i$,
say $F_{i_1},\cdots,F_{i_n}$. 
Then $s\pi_2^{-1}(\tilde{V})=\sum_{k=1}^n\tilde{F}_{i_k}\times_{G_0}G_1%
\times_{G_0}\tilde{V}$, and it suffices to show that each $\tilde{F}_{i_k}\times_{G_0}G_1\times_{G_0}\tilde{V}\rightarrow \tilde{V}$ 
is proper. Now $(s,t)\colon G_1\rightarrow G_0\times G_0$ is proper (Theorem \ref{old}), 
hence so is its pullback 
$(\pi_1,\pi_3)\colon \tilde{F}_{i_k}\times_{G_0}G_1\times_{G_0}\tilde{V}%
\rightarrow F_{i_k}\times\tilde{V}$. Since $\tilde{F}_{i_k}$ is compact, 
the map 
$\tilde{F}_{i_k}\times_{G_0}G_1\times_{G_0}\tilde{V}\rightarrow \tilde{V}$
is also proper, as required.

\paragraph{} Now fix a triangulation ${\cal T}$ adapted to ${\cal U}$, so the 
conditions (\ref{i}) and (\ref{ii}) from the beginning of 
Section \ref{simpset}
are satisfied.  Consider the locally finite cover of $M$ by simplices 
$\sigma\in{\cal T}$ of maximal dimension $n$. Fix for each $\sigma$ a 
chart $U_\sigma\supseteq\sigma$ and a lifting 
$\tilde{\sigma}\subseteq\U_\sigma$
as before. Let $H({\cal T})\subseteq G$ be the full subgroupoid 
constructed from
this cover. Thus $H({\cal T})_0$ is the disjoint sum of the $n$-simplices in 
${\cal T}$, and $H({\cal T})_1$ is constructed as the pullback. 
By Lemma \ref{subeq}, the inclusion $H({\cal T})\subseteq G$ induces an equivalence
of categories of sheaves.

We will show (cf.\ Proposition \ref{H1} below), that the space 
$H({\cal T})_1$
is also a disjoint sum of simplices. To be able to do that, we need a canonical 
way of representing the arrows. We will use a refining atlas and a subdivision
of the triangulation ${\cal T}$ to achieve this.
Let ${\cal V}\succcurlyeq{\cal U}$ be a good atlas refining ${\cal U}$, and 
let ${\cal T}'$ be a subdivision of ${\cal T}$ which is good for ${\cal V}$
(cf. Section \ref{trian}). Thus, each simplex $\tau\in{\cal T}'$ is contained 
in a contractible chart $V_\tau$ in ${\cal V}$. Furthermore, since ${\cal V}$ 
is a good atlas, these charts $V_\tau$ can be chosen in such a way, that 
$V_{\tau'}\subseteq V_\tau$ whenever $\tau'\subseteq \tau$. Moreover, we may 
assume $V_\tau\subseteq U_\sigma$ whenever the simplex $\tau\in{\cal T}'$ is 
contained in the $n$-simplex $\sigma\in{\cal T}$. 
For each $\tau\in{\cal T}'$, choose one $n$-simplex $\sigma\in{\cal T}$ 
containing $\tau$ and an embedding 
$\lambda_{\sigma,\tau}\colon\tilde{V}_\tau\rightarrow\tilde{U}_\sigma$ with
the lifting of $\tau$ in $\tilde\sigma$ in its image. Then let $\tilde\tau=\lambda_{\sigma,\tau}^{-1}(\tau\subseteq\tilde\sigma)$, be the 
lifting of $\tau$ in $\tilde{V}_\tau$.
Next, for every other $n$-simplex $\sigma\in{\cal T}$ containing $\tau$,
choose an embedding 
\begin{equation}\label{lambdasigmatau}
\lambda_{\sigma,\tau}\colon\tilde{V}_\tau\rightarrow\tilde{U}_\sigma%
\mbox{ with }\lambda_{\sigma,\tau}(\tilde\tau)\subseteq \tilde{\sigma}.
\end{equation}
We need some explicit notation for these charts, and denote them by
$$
(\U_\sigma,G_\sigma,\varphi_\sigma\colon\U_\sigma\rightarrow U_\sigma)%
\mbox{ and }(\tilde V_\tau,H_\tau,\psi_\tau\colon\tilde V_\tau\rightarrow V_\tau).
$$
Furthermore, we denote by $H_{\inter{\tau}}\subseteq H_\tau$ the isotropy 
group of the interior of $\tau$ (or of $\tilde{\tau}$).

Note that the space of arrows $H({\cal T})_1$ is the disjoint sum of spaces
$$
H({\cal T})(\sigma_0,\sigma_1)=\{g\in G_1|\, s(g)\in\tilde\sigma_1\subseteq%
\U_{\sigma_1}\mbox{ and }t(g)\in\tilde\sigma_0\subseteq\U_{\sigma_0}\}
$$
of arrows  from $\tilde\sigma_1$ to $\tilde\sigma_0$. Consider one of these spaces.

By assumption on the triangulation ${\cal T}$, there is a family of faces
of $\rho=\sigma_0\cap\sigma_1$,
\begin{equation}\label{filt}
\rho=\rho_k\supseteq\rho_{k-1}\supseteq\cdots\supseteq\rho_0,
\end{equation}
such that the isotropy is constant on $\rho_i-\rho_{i-1}$.
We may assume that $\rho_0=v(\rho)$ is the chosen vertex with maximal isotropy.
By working with the liftings $\tilde\sigma_i\,(i=0,1)$ and the associated liftings
$
\tilde\rho^i=\tilde\rho_k^i\supseteq\tilde\rho_{k-1}^i\supseteq\cdots\supseteq\tilde\rho_0^i
$
contained in $\tilde\sigma_i$, the isotropy subgroups form an increasing family of 
subgroups of the groups $G_{\sigma_i}$,
\begin{equation}\label{subgroups}
G_k^i\subseteq G_{k-1}^i\subseteq\cdots\subseteq G_0^i\subseteq G_{\sigma_i},
\end{equation}
where $G_0^i=G_{v(\rho^i)}$ is the isotropy group of the lifted vertex $\tilde\rho^i_0$ 
of $\tilde\sigma_i$. The filtration (\ref{filt}) yields a similar filtration of 
the space of arrows $H({\cal T})(\sigma_0,\sigma_1)$.
To see this, consider any arrow $g\colon\tilde x_1\rightarrow\tilde x_0$
in $H({\cal T})(\sigma_0,\sigma_1)$. Then in the underlying space $M$
we have $x_0=x_1\in\rho$, so $x_0=x_1$ is contained in one of the simplices
$\tau$ of ${\cal T}'$. It follows that every arrow can be represented in the form 
\begin{equation}\label{arrow}
\U_{\sigma_0}\stackrel{\lambda_0}{\longleftarrow}\tilde V_\tau\stackrel{\lambda_1}{\longrightarrow}\U_{\sigma_1}\quad g=[\lambda_0%
\circ \ell ,\tilde z,\lambda_1]
\end{equation}
where $\lambda_i=\lambda_{\sigma_i,\tau}\,\,(i=0,1)$ are the chosen 
embeddings in (\ref{lambdasigmatau}), $\tilde z\in\tilde\tau\subseteq\tilde V_\tau$ 
is the unique point with $\lambda_i(\tilde{z})=\tilde x_i$, and 
$l\in H_{\tilde{z}}\subseteq H_\tau$ fixes the interior of 
$\tilde\tau$, i.e.\ $l\in H_{\inter{\tau}}$. Let us say that an arrow $g$ 
represented as in (\ref{arrow})
has {\it rank} $\geq k$ if $l\in H_{\inter{\tau}}$ where $\tau\subseteq\rho_{k}$.

For an arrow of rank at least $k$, the following lemma gives a criterion as to
whether the rank is strictly larger than $k$.

\begin{lemma}\label{char}
Consider the open sets
$$
\tilde{N}_i=\lambda_i^{-1}(\tilde{\rho}_{k+1}^i-\tilde\rho_{k}^i)\qquad(i=0,1)
$$
where $\tilde{\rho}_{k+1}^i,\tilde\rho_{k}^i\subseteq\tilde\sigma_i$. An arrow 
$g$ as in (\ref{arrow}) of rank $\geq k$ is of rank $\geq k+1$ if and only if
there is an open neighborhood $W_\tau$ of $\tilde\tau$ in 
$\tilde N_1\cup\tilde\tau$
such that for every $\tilde{y}_1\in (W_\tau-\tau)$, the image 
$\tilde{y}_0=\ell \cdot \tilde{y}_1$ belongs to $\tilde N_0 $.
$$
\makebox[0pt]{\epsfysize=.15\textheight\epsffile{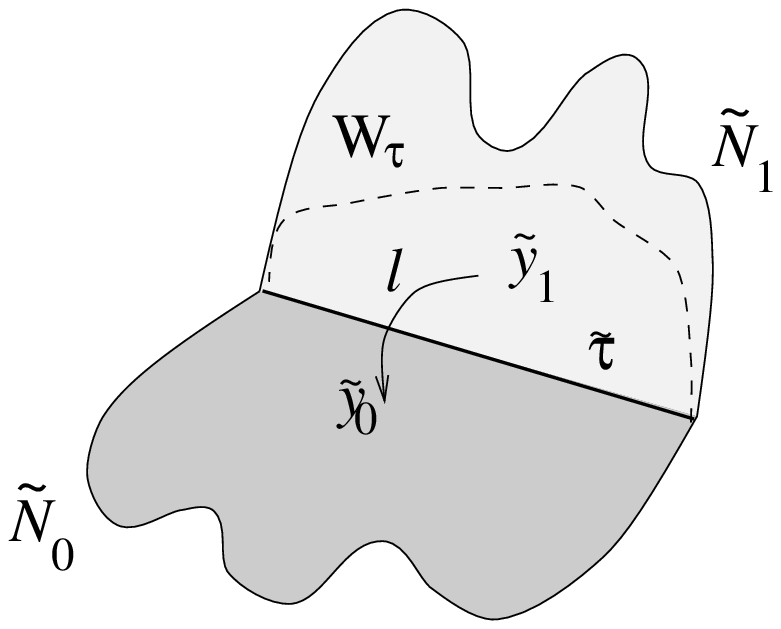} }
$$
\end{lemma}

\paragraph{Proof.}
Suppose that $g$ as represented in (\ref{arrow}) has rank $\geq k+1$.
Let $\tau'\supseteq\tau$ be any simplex in ${\cal T}'$ with 
$\tau'\subseteq\rho_{k+1}$, and consider the following diagram
\begin{equation}\label{diagram}
\diagram
&&\tilde V_\tau\dllto|<\hole|<<\ahook_{\lambda_0}\drrto|<\hole|<<\ahook^{\lambda_1}%
\ddto|<\hole|<<\ahook^\nu\\
\U_{\sigma_0}&a_0&&a_1&\U_{\sigma_1}\\
&&\tilde{V}_{\tau'}\ullto|<\hole|<<\ahook^{\overline{\lambda}_0}%
\urrto|<\hole|<<\ahook_{\overline{\lambda}_1}
\enddiagram
\end{equation}
where $\lambda_i=\lambda_{\sigma_i,\tau}$ as above and 
$\overline{\lambda}_i=\lambda_{\sigma_i,\tau'}$, while $\nu$ 
is any embedding
mapping the chosen lifting $\tilde\tau\subseteq\tilde V_\tau$ into 
$\tilde{\tau}'$. Furthermore,  $a_0,a_1\in H_\tau$ are group elements 
such that
\begin{equation}\label{commelm}
\lambda_i\circ a_i=\overline{\lambda}_i\circ\nu.
\end{equation}
Then by the equivalence relation defining $G_1$, the arrow 
$g=[\lambda_0\circ \ell ,\tilde z,\lambda_1]$ can be represented as
\begin{equation}\label{presentation}
g=[\overline{\lambda}_0\circ \ell',\nu(\tilde z),%
\overline{\lambda}_1]\quad \ell'=\nu(a_0^{-1}\ell a_1),
\end{equation}
and $\ell'\in H_{\interaccent{\tau}}$ by assumption. 
Let $y\in(\tau'-\tau)\cap U_{\sigma_0}\cap U_{\sigma_1}\cap V_{\tau}%
\subseteq(\rho_{k+1}-\rho_k)$ be any point,
and denote its liftings in $\tilde\sigma_i\,(i=0,1)$ by $y_i$.
Let $\tilde{y}_i=\lambda_i^{-1}(y_i)\in\tilde N_i\subset\tilde{V}_\tau$.  
We claim that
$$\ell \cdot \tilde{y}_1 = \tilde{y}_0\;.$$
Indeed, 
$\bar{\lambda}_i \nu (a^{-1}_{i} \tilde{y}_{i}) =
\lambda_{i}(\tilde{y}_i) \in \tilde{\sigma}_i$, and hence
$\nu(a^{-1}_{i}\;\tilde{y}_i) \in \tilde{\tau}'$ because $\bar{\lambda}_i$ is
chosen to map $\tilde{\tau}'$ into $\tilde{\sigma}_i$. It then follows first
that $\nu (a^{-1}_{0}\;\tilde{y}_0) = \nu(a^{-1}_{1}\;\tilde{y}_1)$, since both
sides are liftings in $\tilde{\tau}'$ of the same point $y$, and next that this
point in $\tilde{\tau} '$ is fixed by $\ell'$ since $\ell' \in
H_{\interaccent{\tau}}$.
Thus $\ell' \cdot \nu(a^{-1}_{1} \tilde{y}_{1}) = \nu(a^{-1}_{0}
\tilde{y}_{0})$; or, by definition of $\ell'$ in (15) above, $\nu(a^{-1}_{0}
\ell \tilde{y}_{1}) = \nu(a^{-1}_{0} \tilde{y}_0)$. Since $\nu$  is an
embedding, we conclude that $\ell \tilde{y}_{1} = \tilde{y}_{0}$, as claimed.
So for every $\tau'\supseteq \tau$ with $\tau'\subseteq\rho_{k+1}$ and every $y\in(\tau'-\tau)\cap\varphi_\tau(\tilde N_0)\cap\varphi_\tau(\tilde N_1)$ with 
liftings $\tilde y_i\in\tilde N_i$, we find that 
$\ell\cdot\tilde y_1=\tilde y_0$.
Therefore 
$$
\varphi_\tau^{-1}\left\{\coprod_{\tau\subseteq\tau'\subseteq%
\rho_{k+1}}(\tau'\cap \varphi_\tau(\tilde N_0))\right\}\cap\tilde N_1
$$ 
satisfies the requirements for $W_\tau$. So we have shown that $g$ satisfies 
the condition formulated in this lemma.

Now assume that the arrow $g$ is of rank $\geq k$ and satisfies this condition.
Let $\tau'\supseteq \tau$ be a simplex in $\rho_{k+1}$. Then its inverse 
image in $\tilde N_1$
has a nonempty intersection with $W_\tau$, so its inverse images
in $\tilde N_i$ contain points $\tilde{y}_i$ as in the lemma. 
Use diagram (\ref{diagram}) again to label all the embeddings and 
group elements
involved. 
Then $g$ can again be represented as in (\ref{presentation}) and we need to 
show that $\ell'$ in that presentation is an element of $H_{\interaccent{\tau}}$. 
Let $\tilde{y}_0$ and $\tilde{y}_1$ be as in the lemma, such that $\tilde{y}_i\in\psi_\tau^{-1}(\interaccent{\tau})$.
Then $\ell\cdot\tilde{y}_1=\tilde{y}_0$, so 
$\nu(\ell a_1)\cdot(\overline{\lambda}_1^{\,-1}(\lambda_1(\tilde y_1)))=%
\nu(a_0^{-1})\cdot(\overline{\lambda}_0^{\,-1}(\lambda_0(\tilde y_0)))%
\in\tilde{V}_{\tau'}$, or $\nu(a_0^{-1}\ell a_1)\cdot(\overline{%
\lambda}_1^{\,-1}(\lambda_1(\tilde y_1)))=\overline{\lambda}_0^{\,-1}%
(\lambda_0(\tilde y_0))$.
However, both $\overline{\lambda}_1^{\,-1}(\lambda_1(\tilde y_1))$ and $\overline{\lambda}_0^{\,-1}(\lambda_0(\tilde y_0))$ are liftings in 
$\tilde{\tau}'$ of the same point in $M$, so they have to be the same
in $\tilde{\tau}'$ as well. We conclude that 
$\ell'=\nu(a_0^{-1}\ell a_1)\in H_{\overline{\lambda}_1^{\,-1}(\lambda_1%
(\tilde y_1))}=H_{\interaccent{\tau}}$
(since the isotropy is constant on the interior of a simplex) as required.

\paragraph{}Now we are ready to prove:

\begin{prop}\label{H1}
The space $H({\cal T})_1$ is (homeomorphic to) a disjoint sum of simplices.
\end{prop}

\paragraph{Proof.}
Let $\sigma_0$ and $\sigma_1$ be $n$-simplices in ${\cal T}$ and suppose 
that
$\sigma_0\cap\sigma_1\neq\emptyset$. As we remarked before, it is sufficient 
to prove that $H({\cal T})(\sigma_0,\sigma_1)$ consists of a disjoint sum of 
simplices. Let $\rho_j$ be a part of the filtration (\ref{filt}) above.
Write
\begin{equation}\label{simplices}
\rho_j=\tau_1\cup\cdots\cup\tau_m,
\end{equation}
where $\tau_i$ are simplices in ${\cal T}'$ of the same dimension as $\rho_j$.
Consider all arrows of the form (\ref{arrow}),
$$
\U_{\sigma_0}\stackrel{\lambda_{0,i}}{\longleftarrow}\tilde{V}_{\tau_i}\stackrel%
{\lambda_{1,i}}{\longrightarrow}\U_{\sigma_1}\quad g=[\lambda_{0,i}\circ l,%
\tilde z,\lambda_{1,i}]
$$
where $\lambda_{\varepsilon,i}=\lambda_{\sigma_\varepsilon,\tau_i}$ 
$\varepsilon=0,1;\,(i=1,\cdots,m)$ are the chosen embeddings as in 
(\ref{lambdasigmatau}), and $g$ has rank exactly $j$. For fixed $i$ and 
$\ell $, these arrows form a copy $\tau_i(\ell )$ of the simplex $\tau_i$. 
Moreover, if $\theta=\tau_i\cap\tau_{i'}$ is a nonempty face not in the 
boundary of $\rho_j$ then this copy $\tau_i(\ell )$ is glued (along $\theta$) 
to exactly one copy $\tau_{i'}(\ell')$ in the space 
$H({\cal T})(\sigma_0,\sigma_1)$, as follows.

Since ${\cal T}'$ is a good triangulation, there are embeddings 
$\nu\colon\tilde{V}_\theta\rightarrow\tilde V_{\tau_i}$ and 
$\nu'\colon\tilde V_\theta\rightarrow\tilde V_{\tau_{i'}}$, 
mapping $\tilde\theta$ to $\tilde\tau_i$ and $\tilde\tau_{i'}$, respectively. 
Thus there are $a_0,a_1,b_0,b_1\in H_\theta$ such that
\begin{eqnarray*}
\lambda_{\sigma_j,\theta}a_j&=&\lambda_{j,i}\nu\quad(j=0,1)\\
\lambda_{\sigma_j,\theta}b_j&=&\lambda_{j,i'}\nu'
\end{eqnarray*}
Let $\ell'$ and $\ell''$ be such that $\ell =\nu(\ell'')$ and 
$\ell'=\nu'(b_0^{-1}a_0\ell''a_1^{-1}b_1)$. Then $\tau_i(\ell )$ 
is glued to $\tau_{i'}(\ell')$. (Notice that 
$g=[\lambda_{0,i}\circ \ell ,\tilde z,\lambda_{1,i}]$
has rank exactly $j$ if and only if 
$g'=[\lambda_{0,i'}\circ \ell',\tilde z',\lambda_{1,i'}]$ does. This follows 
from Lemma~\ref{char}, since every pair of open neighborhoods of 
$\tau_i$ and $\tau_{i'}$ have a non-empty intersection. )

Thus, the subspace of all $H({\cal T})(\sigma_0,\sigma_1)$ of all 
these copies $\tau_i(\ell )$ is a covering projection of $\rho_j$, 
hence a disjoint sum of copies of $\rho_j$.

Finally every arrow $g\in H({\cal T})(\sigma_0,\sigma_1)$ occurs in this way,
i.e.\ is represented in the form (\ref{arrow}) where $\tau_i\subseteq\rho_j$
is one of the simplices in (\ref{simplices}) and $g$ has rank exactly $j$.
(This follows easily from considerations as in the proof of Lemma \ref{char}.)

\begin{lemma}
Each space $H({\cal T})_n$ in the nerve of the groupoid $H({\cal T})$ 
is a disjoint sum of simplices.
\end{lemma}

\paragraph{Proof}
This is clear from the fact that $H({\cal T})_0$ and $H({\cal T})_1$ 
are sums of simplices, while $H({\cal T})_n$ is constructed as an 
iterated fibered product
along the source and target maps, which are embeddings on every 
component of $H({\cal T})_1$. 

\paragraph{} Now consider the spectral sequence of Corollary \ref{subspec}.
Let ${\cal A}$ be any locally constant sheaf on ${\cal M}$, and 
let $A$ be the associated $H$-sheaf. Each sheaf $A^{(p)}$ on $H_p$ 
is again locally constant, hence constant on each connected component. 
So in fact $A^{(\bullet)}$ corresponds to a local system of coefficients 
on the simplicial set $\pi_0(H({\cal T})_\bullet)$, obtained by taking the 
connected components of the space $H({\cal T})_n$ in the nerve. Since each 
such 
connected component is a simplex, the spectral sequence of 
Corollary \ref{subspec} collapses, to give the following isomorphism.

\begin{lemma}
For any locally constant sheaf ${\cal A}$ on ${\cal M}$ there is a 
natural isomorphism
$$
H^p(\pi_0(H({\cal T})_\bullet), A)=H^p({\cal M,A}),
$$
where $A$ is the local system of coefficients on the simplicial set 
$\pi_0(H({\cal T})_\bullet)$ induced by ${\cal A}$.
\end{lemma}

The proof of Theorem \ref{main} is now completed by the observation 
that this simplicial set is exactly the simplicial set $S$ described 
in Section \ref{simpcom}.

\begin{lemma}
There is a natural isomorphism of simplicial sets
$$
\pi_0(H({\cal T})_\bullet)\stackrel{\sim}{\longrightarrow}S.
$$
\end{lemma}

\paragraph{Proof.}
By definition, $\pi_0(H({\cal T})_0)=S_0$ in case $n=0$.
For $n>0$ the identity follows from the fact that the definition
of $\mu$ is related to the equivalence relation on 
$H({\cal T})_1\subseteq G_1$ in the following way. Let $\sigma_0$ and 
$\sigma_1$ be $n$-simplices
as before and let $v$ and $w$ be vertices in $\sigma_0\cap\sigma_1$,
connected by a 1-simplex $\theta$ and suppose that $w=v(\theta)$. Then 
$[\lambda_{\sigma_0,v}\circ g,\tilde v,\lambda_{\sigma_1,v}]$ and $[\lambda_{\sigma_0,w}\circ h,\tilde w,\lambda_{\sigma_1,w}]$ are in the same connected component of 
$H({\cal T})_1$ if an only if $h=\mu(g)$. 

\paragraph{} As said, this completes the proof of Theorem \ref{main}.

\vspace*{1cm}

\begin{tabular}{ll}
I. Moerdijk&  D.A. Pronk\\
Mathematical Institute& Department of Mathematics\\
Utrecht University& Dalhousie University\\
PO Box 80010&Halifax, NS\\
3508 TA, Utrecht&Canada, B3H 3J5\\
The Netherlands\\
moerdijk@math.ruu.nl&pronk@cs.dal.ca
\end{tabular}
\end{document}